\documentclass[12pt,reqno]{amsart}
\usepackage{amssymb}

\usepackage[latin1]{inputenc}
\usepackage{mathptmx}

\usepackage{pstricks, pst-3d}
\newcommand{\mm}{\mathfrak m}

\newcommand{\pp}{\mathfrak p}

\newcommand{\Z}{\mathbb{Z}}
\newcommand{\R}{\mathbb{R}}
\newcommand{\N}{\mathbb{N}}
\newcommand{\Q}{\mathbb{Q}}

\newcommand{\Cc}{\mathcal{C}}
\newcommand{\Fc}{\mathcal{F}}

\newcommand{\Kc}{\mathcal{K}}
\newcommand{\Mcc}{\mathcal{M}}
\newcommand{\Ncc}{\mathcal{N}}

\DeclareMathOperator{\pnt}{\raise 0.5mm \hbox{\large\bf.}}
\DeclareMathOperator{\lpnt}{\hbox{\large\bf.}}

\DeclareMathOperator{\gp}{gp}
\DeclareMathOperator{\conv}{conv}

\DeclareMathOperator{\Hilb}{Hilb}

\DeclareMathOperator{\Tor}{Tor}
\DeclareMathOperator{\Ker}{Ker}
\DeclareMathOperator{\projdim}{proj\,dim}
\DeclareMathOperator{\ini}{in}

\DeclareMathOperator{\htt}{ht}

\DeclareMathOperator{\Rad}{Rad}

\def\+#1{\relax\ifmmode\if\noexpand #1\relax \mathop{\kern
    0pt^+{#1}}\nolimits\else \kern 0pt^+\!#1 \fi\else$^*$#1\fi}

\def\ie{i.\,e.\ }
\def\fg{f.\,g.\ }
\def\eg{e.\,g.\ }
\def\Eg{E.\,g.\ }

\let\phi=\varphi
\let\:=\colon

\newtheorem{thm}{\bf Theorem}[section]
\newtheorem{lem}[thm]{\bf Lemma}
\newtheorem{cor}[thm]{\bf Corollary}
\newtheorem{prop}[thm]{\bf Proposition}

\theoremstyle{definition}

\newtheorem{rem}[thm]{\bf Remark}
\newtheorem{ex}[thm]{\bf Example}

\theoremstyle{plain}
\newtheorem*{thm*}{Theorem}

\title{Gröbner bases and Betti numbers of monoidal complexes}

\author{Winfried Bruns}
\address{Fachbereich Mathematik/Informatik, Institut für Mathematik, Universität Osna\-brück, Albrechtstr. 28a, 49069 Osnabrück, Germany}
\email{wbruns@uos.de}

\author{Robert Koch}\address{An der Kommende 11, 53177 Bonn, Germany} \email{rrkoch@web.de}

\author{Tim Römer}
\address{Fachbereich Mathematik/Informatik, Institut für Mathematik, Universität Osna\-brück, Albrechtstr. 28a, 49069 Osnabrück, Germany}
\email{troemer@uos.de}

\dedicatory{To Mel Hochster on his 65th birthday}

\begin{document}

\begin{abstract}
In this note we consider monoidal complexes and their associated
algebras, called toric face rings. These rings generalize
Stanley-Reisner rings and affine monoid algebras. We compute initial
ideals of the presentation ideal of a toric face ring, and determine
its graded Betti numbers. Our results generalize celebrated theorems
of Hochster in combinatorial commutative algebra.
\end{abstract}

\maketitle
%
%
\section{Introduction}
\label{intro} Combinatorial commutative algebra is a branch of
combinatorics, discrete geometry and commutative algebra. On the one
hand problems from combinatorics or discrete geometry are studied
using techniques from commutative algebra. On the other hand
questions in combinatorics motivated various results in commutative
algebra. Since the fundamental papers of Stanley (see \cite{ST05}
for the results) and Hochster \cite{HO72, HO77} combinatorial
commutative algebra is a growing and very active field of research.
See also Bruns-Herzog \cite{BRHE98}, Villarreal \cite{VI1},
Miller-Sturmfels \cite{MIST05} and Sturmfels \cite{ST96} for
classical and recent results and new developments in this area of
mathematics.

Stanley-Reisner rings and affine monoid algebras are two of the
classes of rings considered in combinatorial commutative algebra. In
this paper we consider toric face rings associated to monoidal
complexes. They generalize Stanley-Reisner rings by allowing a more
general incidence structure than simplicial complexes, and
more general rings associated with their faces, namely affine monoid
algebras instead of polynomial rings.

In cooperation with M. Brun and B. Ichim the authors have studied
the local cohomology of toric face rings in previous work
\cite{BRBRRO07}, \cite{BRRO08}, \cite{ICRO}, and one of the main
results is a general version of Hochster's formula for the local
cohomology of a Stanley-Reisner ring (see \cite{BRHE98} or
\cite{ST05}), even beyond toric face rings.

In this paper we want to generalize Hochster's formulas for the
graded Betti numbers of a Stanley-Reisner ring \cite{HO77} and
affine monoid rings \cite[Theorem 9.2]{MIST05} to toric face rings. Such a
generalization is indeed possible for monoidal complexes, that,
roughly speaking, can be embedded into a space $\Q^d$. As
counterexamples show, full generality does not seem possible. One of
the problems encountered is to construct a suitable grading. This
forces us to consider grading monoids that are not necessarily
cancellative.

Another topic treated are initial ideals (of the defining ideals) of
toric face rings with respect to monomial (pre)orders defined by
weights. Indeed, toric face rings come up naturally in the study of
initial ideals of affine monoid algebras. In this
regard we generalize results of Sturmfels \cite{ST96}. We will pay
special attention to the question when the initial ideal is radical,
monomial or both. This gives an opportunity to indicate a
``simplicial'' proof of Hochster's famous theorem on the
Cohen-Macaulay property of affine normal monoid domains \cite{HO72}.
For unexplained terminology we refer the reader to \cite{BRGU08} and
\cite{BRHE98}.


\section{Monoidal complexes and toric face rings}
\label{facerings}

A \emph{cone} is a subset of a space $\R^d$ of type
$\R_+x_1+\dots+\R_+x_n$ with $x_1,\dots,x_n\in\R^d$. The
\emph{dimension} of a cone $C$ is the vector space dimension of $\R
C$. A \emph{face} of $C$ is a subset of type $C\cap H$ where $H$ is
a \emph{support hyperplane} of $C$, \ie a hyperplane $H$ for which
$C$ is contained in one of the two closed halfspaces $H^+, H^-$
determined by $H$. A \emph{rational} cone is generated by elements
$x\in\Q^d$. A {\em pointed} cone has $\{0\}$ as a face.

A \emph{fan} in $\R^d$ is a finite collection $\Fc$ of cones in
$\R^d$ satisfying the following conditions:
\begin{enumerate}
\item
all the faces of each cone $C\in\Fc$ belong to $\Fc$, too;
\item
the intersection $C\cap D$ of $C,D\in\Fc$ is a face of $C$ and of
$D$.
\end{enumerate}

We want to investigate more general configurations of cones, giving
up the condition that all cones are contained in a single space, but
retaining the incidence structure. A {\em conical complex} consists
of
\begin{enumerate}
\item
a finite set $\Sigma$ of sets,
\item
a cone $C_c \subseteq \R^{\delta_c}$, $\delta_c=\dim \R C_c$, for
each $c \in \Sigma$,
\item
and a bijection $\pi_c \: C_c \to c$ for each $c \in \Sigma$ such
that the following conditions are satisfied:
\begin{enumerate}
\item
for each face $C'$ of $C_c$, $c \in \Sigma$, there exists $c' \in
\Sigma$ with $\pi_c(C')=c'$;
\item
for all $c,d \in \Sigma$ there exist faces $C'$ of $C_c$ and $D'$ of
$D_d$ such that $c\cap d=\pi_c(C') \cap \pi_d(D')$ and the
restriction of $\pi_d^{-1}\circ \pi_c$ to $C'$ is an isomorphism of
the cones $C'$ and $D'$.
\end{enumerate}
\end{enumerate}
Here
an \emph{isomorphism} of cones $C,D$ is a bijective map $\phi\colon C\to D$
that extends to an isomorphism of the vector spaces $\R C$ and $\R
D$. Simplifying the notation, we write $\Sigma$ also for the conical
complex. A fan $\Fc$ is a conical complex in a natural way: fans are
nothing but embedded conical complexes.

As introduced in the definition, $\delta_c$ will always denote the
dimension of $C_c$ so that $\R C_c$ can be identified with
$\R^{\delta_c}$. The elements $c\in \Sigma$ are called the {\em
faces} of $\Sigma$. Similarly, one defines {\em rays} and {\em
facets} of $\Sigma$ as $1$-dimensional and maximal faces of
$\Sigma$. The {\em dimension} of $\Sigma$ is the maximal dimension
of a facet of $\Sigma$. We denote by $|\Sigma|=\bigcup_{c \in
\Sigma} c$ the {\em support} of $\Sigma$. Identifying $C_c$ with
$c$, we may consider $C_c$ as a subset of $|\Sigma|$. Then we can
treat $|\Sigma|$ almost like an (embedded) fan. The main difference
is that it makes no sense to speak of concepts like convexity
globally. However, locally in the cones $C_c$ we may consider convex
subsets. The complex $\Sigma$ is {\em rational} and {\em pointed},
respectively, if all cones $C_c$, $c \in \Sigma$, are rational and
pointed respectively. We call $\Sigma$ {\em simplicial}, if all
cones $C_c$, $c \in \Sigma$, are simplicial, \ie they are generated
by linearly independent vectors.

In order to define interesting algebraic objects associated to a
conical complex one needs a corresponding discrete structure. A {\em
monoidal complex} $\Mcc$ supported by a conical complex $\Sigma$ is
a set of monoids $(M_c)_{c \in \Sigma}$ such that
\begin{enumerate}
\item
for each $c\in C$ the monoid $M_c$ is an affine (\ie finitely
generated) monoid contained in $\Z^\delta_c$;

\item $M_c
\subseteq C_c$ and $\R_+ M_c=C_c$ for every $c \in \Sigma$;

\item
for all  $c,d \in \Sigma$ the map $\pi_d^{-1}\circ \pi_c$ restricts
to a monoid isomorphism between $M_c \cap \pi_c^{-1}(c\cap d)$ and
$M_d \cap \pi_d^{-1}(c\cap d)$.
\end{enumerate}
In other words, we have chosen an affine monoid $M_c$ for every $c
\in \Sigma$ which generates $C_c$ and whose intersection with a face
$C_d$ of $C_c$ is just $M_d$. The monoidal complex naturally
associated to a single affine monoid $M$ is simply denoted by $M$;
it is supported on the conical complex formed by the faces of the
cone $\R_+ M$.

The simplest examples of \emph{conical complexes} are those
\emph{associated with rational fans $\Fc$}. For each cone $C\in\Fc$
we choose $M_C=C\cap \Z^d$. These monoids are finitely generated by
Gordan's lemma. Moreover, they are normal: recall that an affine
monoid $M$ is \emph{normal} if $M=\gp(M)\cap \R_+M$.

\begin{rem}
Let $\gp(M)$ denote the group of differences of a monoid $M$. The
groups $\gp(M_C)$ of the monoids in a monoidal complex associated
with a fan form again a
monoidal complex in a
 natural way since $\gp(M_D)=\gp(M_C)\cap \R
D$ if $D$ is a face of $C$.

In general the compatibility condition between the passage to faces
and the formation of groups of differences need not be satisfied.
Nevertheless, the rational structures defined by the monoids $M_c$,
namely the rational subspaces $\Q \gp(M_c)$ of $\R^{\delta_c}$ are
compatible with the passage to faces. This follows from condition
(ii): both monoids $\gp(M_c)\cap \R C_d$ and $\gp(M_d)$ are
contained in $\Z^{\delta_d}$ and have the same rank $\delta_d$.
\end{rem}

Note that the monoids $M_c$ form a direct system of sets with
respect to the embeddings $\pi_d^{-1}\circ \pi_c \: M_c \to M_d$
where $c,d \in \Sigma$ and $c \subseteq d$. We set
$$
|\Mcc| = \varinjlim \ M_c.
$$
In general, there exists no global monoid structure on $|\Mcc|$, but
it carries a partial monoid structure since we can consider each
monoid $M_c$ as a subset of $|\Mcc|$ in the natural way. Whenever
there exists $c\in \Sigma$ such that $a,b\in M_c$ then $a+b$ is
their sum in $M_c$, and as an element of $|\Mcc|$ the sum is
independent of the choice of $c$.

Next we choose a field $K$ and define the {\em toric face ring}
$K[\Mcc]$ of $\Mcc$ (over $K$) as follows. As a $K$-vector space let
$$
K[\Mcc] = \bigoplus_{a \in |\Mcc|} K t^a.
$$
We set
$$
t^a \cdot t^b =
\begin{cases}
t^{a+b} & \text{if $a, b \in M_c$ for some $c \in \Sigma$},\\
0 & \text{otherwise}.
\end{cases}
$$
Multiplication in $K[\Mcc]$ is defined as the $K$-bilinear extension
of this product. It turns $K[\Mcc]$ into a $K$-algebra. In the
following, the elements of $|\Mcc|$ are called {\em monomials}.

There exist at least two other natural descriptions of toric face
rings of a monoidal complex. The first is a realization as an
inverse limit of the affine monoid rings $K[M_c]$, $c \in \Sigma$.
For $c\in \Sigma$ and a face $d$ of $c$ there exists a natural
projection map $K[M_c] \to K[M_{d}]$ which sends monomials $t^a$ to
zero if $a \notin M_{d}$, the {\em face projection map}. With
respect to these maps we may consider the inverse limit $\varprojlim
\ K[M_c]$:

\begin{prop}
\label{inverselimit} Let $\Mcc$ be a monoidal complex supported on a
conical complex $\Sigma$. Then
$$
K[\Mcc] \cong \varprojlim \ K[M_c].
$$
\end{prop}

For the proof of the proposition we introduce some more notation.
Let $c \in \Sigma$ and let $\pp_c$ be the ideal of $K[\Mcc]$ which
is generated by all monomials $t^a$ with $a\notin M_c$. Then there
is a natural isomorphism of $K$-algebras $K[M_c] \cong
K[\Mcc]/\pp_c$. In particular, $\pp_c$ is a prime ideal. Moreover,
if $d\subset c$, $c,d \in \Sigma$, then the natural epimorphism
$K[\Mcc]/\pp_c \to K[\Mcc]/\pp_d$ coincides with the map induced
from the projection map $K[M_c] \to K[M_d]$, and we identify these
maps in the following.

\begin{proof}[Proof of Proposition \ref{inverselimit}]
Observe that each of the ideals $\pp_c$ has a $K$-basis consisting
of monomials of $K[\Mcc]$. Therefore the following equations are
satisfied for $c,d,e \in \Sigma$:
\begin{enumerate}
\item
$\pp_c+\pp_d=\pp_{c\cap d}$,
\item
$\pp_c\cap(\pp_d+\pp_e)=\pp_c\cap \pp_d+\pp_c\cap \pp_e$,
\item
$\pp_c+\pp_d\cap \pp_e=\pp_c\cap \pp_d+\pp_c\cap \pp_e$ for all $i,j,k$.
\end{enumerate}
Now it follows easily that $\varprojlim K[\Mcc]/\pp_c$ is isomorphic
to $K[\Mcc]/\bigcap_{c\in\Gamma}\pp_c$ (for example, see Example 3.3
in \cite{BRBRRO07}). But $\bigcap_{c\in\Gamma}\pp_c=0$, and so
\begin{equation*}
 \varprojlim K[M_c] \cong
 \varprojlim K[\Mcc]/\pp_c \cong
 K[\Mcc]/\bigcap_{c\in\Gamma}\pp_c\cong
 K[\Mcc].\qedhere
\end{equation*}
\end{proof}

Second, we want to describe a toric face ring as a quotient of a
polynomial ring. It is not difficult to compute the defining ideal
of such a presentation. In view of Theorem \ref{Sturm1} below we
have to consider elements of $K[\Mcc]$ that are either monomials
$t^a$, $a\in\Mcc$, or $0$. For a uniform notation we augment
$|\Mcc|$ by an element $-\infty$ and set $t^{-\infty}=0$.

\begin{prop}
\label{presentation} Let $\Mcc$ be a monoidal complex supported on a
conical complex  $\Sigma$, and let $(a_e)_{e\in E}$ be a family of
elements of $|\Mcc|\cup\{-\infty\}$ generating $K[\Mcc]$ as an
$K$-algebra. (Equivalently, $\{a_e:e\in E\}\cap M_c$ generates $M_c$
for each $c\in\Sigma$.) Then the kernel $I_{\Mcc}$ of the surjection
$$
\phi\: K[X_e:e\in E] \to K[\Mcc],\qquad \phi(X_e)=t^{a_e},
$$
is generated by
\begin{enumerate}
\item
all monomials $\prod_{h\in H} X_h$ where $H$ is a subset of $E$ for
which $\{a_h : h\in H\}$ is not contained in any monoid $M_c$, $c\in
\Gamma$, and
\item
all binomials $\prod_{g \in G} X_g^{i_g} - \prod_{h \in H}
X_h^{j_h}$ where $G,H\subset E$, all $a_g, a_h$ are contained in a
monoid $M_c$ for some $c\in \Sigma$, and $\sum_{g\in G} i_g
a_g=\sum_{h\in H} j_h a_h$.
\end{enumerate}
Moreover, the monomials in (i) are all monomials contained in
$I_{\Mcc}$. A binomial is contained in $I_{\Mcc}$ if either both its
monomials are contained in the family of monomials given in (i), or
it is in the list of the binomials in (ii).
\end{prop}

\begin{proof}
It is clear that $I_{\Mcc}$ contains the ideal $J$ generated by all
the monomials and binomials listed in (i) and (ii).

For the converse, let $f$ be a polynomial such that $\phi(f)=0$.
Then we can assume that all monomials of $f$ map to elements of
$|\Mcc|$ since all other monomials belong to $I_{\Mcc}$. Now let
$c\in \Sigma$, and define $f_c$ to be the polynomial that arises as
the sum of those terms of $f$ whose monomials are mapped to elements
of $M_c \subset |\Mcc|$. Then $\phi(f_c)=0$ as well. It is
well-known and easy to show that  $f_c$ then belongs to the ideal in
$R[X_e: e\in E]$ generated by all those binomials in (ii) for which
$a_g,a_h \in M_c$. (Equivalently, the binomials in (ii) for $M_c$
generate the presentation ideal of  $K[M_c]$ over a polynomial ring
in the variables $X_e$ where $e\in E$ and  $a_e \in M_c$). Therefore
we may replace $f$ by $f-f_c$, and finish the proof by induction on
the number of terms of $f$.

It is clear that a monomial belongs to $I_{\Mcc}$ if and only if it
is contained in the family of monomials given in (i). If a binomial
is an element of $I_{\Mcc}$, then either both monomials belong to
this ideal, or none of the monomials. In the latter case it must be
one of the binomials of the family of binomials given in (ii), since
no other binomials belong to the kernel of the map $K[X_e:e \in E]
\to K[\Mcc]$.  This follows directly from the construction of the
ring $K[\Mcc]$.
\end{proof}

\begin{ex}
\label{ExMonComp} \
\begin{enumerate}
\item
Let $\Fc$ be a rational fan in $\R^d$, and let $\Mcc$ be the conical
complex associated with it. Then the algebra $K[\Mcc]$ is the toric
face ring introduced by Stanley \cite{ST87}.
\item
Let $\Delta$ be an abstract simplicial complex on the vertex set
$[n]=\{1,\dots,n\}$. Then $\Delta$ has a geometric realization by
considering the simplices $\conv(e_{i_1},\dots,e_{i_m})$ such that
$\{i_1,\dots,i_m\}$ belongs to $\Delta$ (here $e_1,\dots,e_n$ is the
canonical basis of $\R^n$). The cones over the faces of the
geometric realization form a fan $\Fc$, and its toric face ring $R$
given by (i) is nothing but the Stanley-Reisner ring of $\Delta$. In
fact, according to Proposition \ref{presentation}, the kernel of the
natural epimorphism $K[X_1,\dots,X_n]\to R$ is generated by those
monomials $X_{j_1}\cdots X_{j_r}$ such that $\{j_1,\dots,
j_r\}\notin \Delta$.

Algebras associated with monoidal complexes therefore generalize
Stanley-Reisner rings by allowing arbitrary conical complexes as
their combinatorial skeleton and, consequently, monoid algebras as
their ring-theoretic flesh.
\item
The polyhedral algebras of \cite{BRGU01} are another special case of
the algebras associated with monoidal complexes. For them the cones
are generated by lattice polytopes and the monoids are the polytopal
monoids considered on \cite{BRGU01}.
\end{enumerate}
\end{ex}

\begin{rem}
\label{othersources} In \cite{BRRO05} toric face rings were defined
by their presentation ideals given in Proposition
\ref{presentation}. Thus Proposition \ref{inverselimit} is
equivalent to \cite[Theorem 4.7]{BRRO05}. In \cite{BRRO05} the
generators of the affine monoids $M_c$ for $c \in \Sigma$ were fixed
right in the beginning while in this paper we fix only the monoids
and are free to choose generators whenever we like. The new approach
leads directly to the natural description of the toric face rings.
Using arguments as above (\eg the prime ideals as in
\ref{inverselimit}) one obtains alternative and slightly more
compact proofs than those in \cite{BRRO05}.
\end{rem}

We have already used the fact that the zero ideal of $K[\Mcc]$ is
the intersection of the prime ideals $\pp_c$. This implies that
$K[\Mcc]$ is reduced.

Let $\Sigma$ be a conical complex. A conical complex $\Gamma$ is a
{\em subdivision} of $\Sigma$ if $|\Gamma|=|\Sigma|$ and each face
$c\in \Sigma$ is the union of faces $d \in \Gamma$. The subdivision
is called a {\em triangulation} if $\Gamma$ is simplicial. We call a
subdivision $\Gamma$ {\em rational}, if all cones $C_d$, $d\in
\Gamma$ are rational.

Suppose that $\Gamma$ is a subdivision of $\Sigma$, let $\Mcc$ be a
monoidal complex supported by $\Gamma$, and $c$ a face of $\Sigma$.
In the situation of Proposition \ref{presentation} for the toric
face ring $K[\Mcc]$ we let $S_c$ be the polynomial subring of
$S=K[X_e:e\in E]$ generated by those $X_e$ for which $a_e\in C_c$.
Furthermore let $\Mcc_c$ be the monoidal subcomplex of $\Mcc$
consisting of all faces $D_d$ of $\Gamma$, $d\subset c$, and their
associated monoids.

Since $\Mcc_c$ is a monoidal subcomplex, one has a natural
epimorphism $K[\Mcc]\to K[\Mcc_c]$, generalizing the face
projection. It is given by $t^a\mapsto t^a$ whenever $a\in C_c$, and
$t^a\mapsto 0$ otherwise.

But we have also an embedding $K[\Mcc_c]\to K[\Mcc]$, since points
of $|\Mcc_c|$ that are contained in a face of $\Gamma$, are also
contained in a face of $\Mcc_c$.

In order to encode the incidence structure of $\Sigma$ we let
$A_\Sigma$ denote the ideal in $S$ generated by the squarefree
monomials $\prod_{h\in H} X^h$ for which $\{a_h:h\in H\}$ is not
contained in a face of $\Sigma$.

\begin{prop}\label{Mc_retr}
With the notation introduced in Proposition \ref{presentation} we
have:
\begin{enumerate}
\item The embedding $K[\Mcc_c]\to K[\Mcc]$ is a section of the
projection $K[\Mcc]\to K[\Mcc_c]$, and thus makes $K[\Mcc_c]$ a
retract of $K[\Mcc]$.

\item Let $c_1,\dots,c_n$ be the facets of $\Sigma$, and set
$\Mcc_i=\Mcc_{c_i}$. Then
$$
I_\Mcc=A_\Sigma+ S I_{\Mcc_{1}}+\dots +S I_{\Mcc_{n}}.
$$
Moreover, for each face $c\in\Sigma$ we have $I_{\Mcc_c}=S_c\cap
I_\Mcc$.
\end{enumerate}
\end{prop}

\begin{proof}
Part (i) is evident, and the representation of $I_\Mcc$ in part
(ii) follows immediately from Proposition \ref{presentation}: none
of the binomial relations is lost on the right hand side, and it
contains also all the monomial relations because these are either
contained in one of the $I_{\Mcc_i}$ or in $A_\Sigma$. The equation
$I_{M_c}=S_c\cap I_\Mcc$ restates part (i), lifted to the
presentations of the algebras.
\end{proof}

In particular we can apply Proposition \ref{Mc_retr} in the case
$\Gamma=\Sigma$.

\section{Toric face rings and initial ideals}
\label{initialideals}

Next we want to compute initial ideals of the presentation ideals of
monoidal complexes considered in Proposition \ref{presentation}.
Recall that a {\em weight vector} for a polynomial ring
$S=K[X_1,\dots,X_n]$ is an element $w \in \N^n$ where $\N$ denotes
the set of non-negative integers. Given this vector we assign $X_i$
the {\em weight} $w_i$. It is easy to see that this is equivalent to
endow $S$ with a positive $\Z$-grading under which the monomials are
homogeneous. Thus the whole terminology of graded rings (with the
prefix $w$) can be applied. In particular, we can speak of the
{\em $w$-degree} of a monomial; it is defined by
$$
\deg_w X^a=\sum_{i=1}^n a_i w_i = a \cdot w.
$$
A weight vector $w$ determines a {\em weight (pre-)order} if one sets
$$
X^a \le_w X^b \quad\iff\quad a \cdot w \le b \cdot w.
$$
The only axiom of a monomial order (as considered below) not
satisfied is antisymmetry: for $n>1$ there always exist distinct
monomials $X^a$ and $X^b$ such that simultaneously $X^a\le_w X^b$
and $X^b\le_w X^a$.

The {\em $w$-initial component} $\ini_w(f)$ of a polynomial $f$ is
simply its $w$-homogene\-ous component of highest degree. Let $V
\subseteq S$ be a subspace. Then the {\em $w$-initial subspace
$\ini_w(V)$} is the subspace generated by the polynomials
$\ini_w(f)$, $f\in V$. Observe that for an ideal $I \subseteq S$ the
$w$-initial subspace $\ini_w(I)$ is again an ideal of $S$. Now
well-known results for monomial orders (see below) hold also for
weight orders. E.g.\ for subspaces $V_1\subseteq V_2 \subseteq S$ we
have $\ini_w(V_1)=\ini_w(V_2)$ if and only if $V_1=V_2$.

A {\em monomial order} $<$ on $S$ is a total order of the monomials
of $S$ such that $1<X^a$ for all monomials $X^a$, and $X^a<X^b$
implies $X^{a+c}<X^{b+c}$ for all monomials $X^a,X^b,X^c$. Now we
can speak similarly of initial terms $\ini_<(f)$ and initial
subspaces $\ini_<(V)$ with respect to $<$. Recall that a {\em
Gröbner basis} of $I$ is a set of elements of $I$ whose initial
monomials generate $\ini_<(I)$. Such a set always exists and then
also generates $I$.

It is an important fact that a monomial order can always be
approximated by a weight order if only finitely many monomials are
concerned: for an ideal $I$ of $S$ there exists a weight vector $w
\in \N^n$ such that $\ini_<(I)=\ini_w(I)$. Conversely, given a
weight vector $w \in \N^n$ and a monomial order $<'$ we can refine
the weight order $<_w$ to a monomial order $<$ by setting $X^a <
X^b$ if either $a\cdot w < b\cdot w$, or $a \cdot w = b \cdot w$ and
$X^a <' X^b$. Observe also that the $w$-initial terms of a Gröbner
basis of $I$ with respect to $<$ generate $\ini_w(I)$. For more
details and general results on weight orders and monomial orders we
refer to \cite{BRCO03} or \cite{ST96}.

The ideal given in Proposition \ref{presentation} has a special
structure. It is generated by monomials and binomials. This property
persists in the passage to initial ideals.

\begin{lem}
\label{IniGenBin} Let $I \subset K[X_1,\dots,X_n]$ be an ideal
generated by monomials and binomials and $w \in \N^n$ a weight
vector. Then $\ini_w(I)$ is generated by the monomials and the
initial components of the binomials in $I$.
\end{lem}
\begin{proof}
We refine the weight order to a monomial order $<$. Using the
Buchberger algorithm to compute a Gröbner basis for $I$, one
enlarges the given set of generators of $I$ consisting of monomials
and binomials only by more monomials and binomials. The
corresponding initial components with respect to the weight order
$<_w$ then generate $\ini_w(I)$.
\end{proof}

It is a useful consequence of Lemma \ref{IniGenBin} that the
decomposition of the ideal $I_\Mcc$ in Proposition \ref{Mc_retr} is
passed onto  their initial ideals. We need it only for the
trivial subdivision of $\Sigma$ by itself, but it can easily be
generalized to the setting of Proposition \ref{Mc_retr}.
(Also see \cite[Theorem 5.9]{BRRO08} for a related result.)

\begin{prop}\label{ini_retr}
Consider the presentation of $K[\Mcc]$ as a residue class ring of
$S=K[X_e:e\in E]$ as in Proposition \ref{presentation}, a weight
vector $w$ on $S$ and the induced weight vectors for the subalgebras
$S_c=K[X_e: a_e\in M_c]$, $c\in \Sigma$. Then
$$
\ini_w(I_\Sigma)=A_\Mcc+S \cdot\ini_w(I_{M_1})+\dots+ S
 \cdot\ini_w(I_{M_n})
$$
where again $c_1,\dots,c_n$ are the facets of $\Sigma$, and
$M_i=M_{c_i}$. Moreover, $\ini_w(I_{M_c})=S_c\cap \ini_w(I_\Mcc)$
for all $c\in\Sigma$.
\end{prop}

\begin{proof}
It is clear that the right hand side is contained in
$\ini_w(I_\Sigma)$. For the converse inclusion it is enough to
consider the system of generators of $I_\Mcc$ described in
Proposition \ref{IniGenBin}, and there is nothing to say about the
monomials in $I_\Mcc$. Let $f$ be the initial component of a
binomial $g$ in $I_\Mcc$. According to Proposition
\ref{presentation} there are two cases: (1) $g$ belongs to $A_\Mcc$;
then so does $f$. (2) $g\in I_{M_i}$ for some $i$; then $f\in
\ini_w(I_{M_i})$, and we are done with the decomposition of
$\ini_w(I_\Sigma)$.

The equality $\ini_w(I_{M_c})=S_c\cap \ini_w(I_\Mcc)$ is left to the
reader. It is easily derived from Proposition \ref{presentation} and
Proposition \ref{IniGenBin}.
\end{proof}

Recall that a function $f \: X \to \R$ on a convex set $X$ is called
{\em convex} if $f(tx+(1-t)y) \leq tf(x)+(1-t)f(y)$ for all $x,y \in
X$ and $t \in [0,1]$. A function $f \: |\Sigma| \to \R$ on a conical
complex $\Sigma$ is called {\em convex} if it is convex on all the
cones $C_c$ for $c \in \Sigma$. For a function $f \: |\Sigma| \to
\R$ a connected subset $W$ of a facet $C_c$ of $\Sigma$ is a {\em
domain of linearity} if it is maximal with respect to the following
property: $g_{|W_c}$ can be extended to an affine function on $\R
C_c$. Now a subdivision $\Gamma$ of a conical complex $\Sigma$ is
said to be {\em regular}, if there exists a convex function $f \:
|\Sigma| \to \R$ whose domains of linearity are facets of $\Gamma$.
Such a function is called a {\em support function} for the
subdivision $\Gamma$.

Let $(a_e)_{e\in E}$ be a family of elements of $|\Mcc|$ such that
$\{a_e:e\in E\}\cap M_c$ generates $M_c$ for each $c\in\Sigma$. Now
we choose a polynomial ring $S=K[X_e: e \in E]$ and define the
surjective homomorphism $\phi\: K[X_e:e\in E] \to K[\Mcc]$ which
maps $X_e$ to $t^{a_e}$ as considered in Proposition
\ref{presentation}. Let $w=(w_e)_{e\in E}$ be a weight vector for
$S$.

On the one hand, the weight vector $w$ determines initial ideals,
especially the initial ideal $\ini_w(I_{\Mcc})$. On the other hand,
$w$ determines also a conical subdivision $\Gamma_w$ of the conical
complex $\Sigma$ as follows. First, every cone $C_c \subseteq
\R^{\delta_c}$ and the weight vector $w$ define the cone
$$
C'_c=\R_+((a_e,w_e): e\in E \text{ such that } a_e \in C_c )
\subseteq \R^{\delta_c+1}.
$$
The projection on the first $\delta_c$ coordinates maps $C'_c$ onto
$C_c$. The \emph{bottom} of $C'_c$ with respect to $C_c$ consists of
all points $(a, h_a)\in C'_c$ such that the line segment
$[\bigl(a,0),(a,h_a)\bigr]$ intersects $C'_c$ only in $(a,h_a)$. In
other words, $h_a=\min\{h': (a,h')\in C'_c\}$. Clearly $h_a>0$ for
all $a\in C_c$, $a \neq0$. The bottom is a subcomplex of the
boundary of $C'_c$ (or $C'_c$ itself). Note that its projection onto
$C_c$ defines a conical subdivision of the cone $C_c$. Second, the
collection of these conical subdivisions of the cones $C_c$
constitutes a conical subdivision of $\Sigma$.

Now we show that this subdivision is regular as defined above. To
this end we define the function $\htt_w \: |\Sigma| \to \R$ as
follows. For $a \in |\Sigma|$ there exists a minimal face $c \in
\Sigma$ such that $a \in C_c$. Construct $C'_c$ as above using the
weight vector $w$. Then we define
$$
\htt_w(a)= \min\{ h' \in \R: (a,h') \in C'_c\},
$$
\ie  $\htt_w(a)$ is the unique vector in the bottom of $C'_c$ which
is projected on $C_c$ via the projection map on the first $\delta_c$
coordinates.

\begin{prop}
\label{htt_conv} Let $\Mcc$ be a monoidal complex supported on a
conical complex  $\Sigma$, let $(a_e)_{e\in E}$ be a family of
elements of $|\Mcc|$ such that $\{a_e:e\in E\}\cap M_c$ generates
$M_c$ for each $c\in\Sigma$, and let $w=(w_e)_{e\in E}$ be a weight
vector. Then:
\begin{enumerate}
\item
For $c \in \Sigma$, $b_1,\dots,b_m \in C_c$ and $\alpha_i>0$,
$i=1,\dots,m$, we have
\begin{equation}
\htt_w\bigl(\sum_{i=1}^m \alpha_i b_i\bigr)\leq \sum_{i=1}^m
\alpha_i \htt_w(b_i).  \label{eq_htt_conv}
\end{equation}
In particular, $\htt_w$ is a convex function on $|\Sigma|$.
\item
Its domains of linearity are the cones $D_d$ for facets $d$ of
$\Gamma_w$, \ie equality holds in \eqref{eq_htt_conv} if and only if
there exists a facet of $\Gamma_w$ containing $b_1,\dots,b_m$.
\end{enumerate}
Therefore $\Gamma_w$ is a regular subdivision of $\Sigma$.
\end{prop}

Part (i) uses only the definition of $\htt$ and that
the cones $C_c'$ are closed under $\R_+$-linear combinations, and
part (ii) reflects the fact that an $\R_+$-linear combination of
points in the boundary of a cone $C$ lies in the boundary if and
only if all points (with nonzero coefficients) belong to a facet of
$C$. (Also see \cite[Lemma 7.16]{BRGU08}.)

Since the weights $w_e$ are positive, the cone $C'_c$ are pointed,
even if $C_c$ is not. Thus all faces of $\Gamma_w$ are pointed, too.

For each $D_d$ with  $d\in \Gamma_w$ we let $N_{d,w}$ be the monoid
generated by all $a_e \in D_d$ for which $\htt_w(a_e) = w_e$. The
cones $D_d$ and the monoids $N_{d,w}$ form a monoidal complex
$\Mcc_{w}=\Mcc_{\Gamma_w}$ supported by the conical complex
$\Gamma_w$, the {\em monoidal complex defined by $w$}. Observe that
each extreme ray of a cone $D_d$ of $\Gamma_w$ is the image of an
extreme ray of $C_c'$ for some $c\in \Sigma$. The latter contains a
point $(a_e,w_e)$, and therefore $w_e=\htt_w(a_e)$. This implies
$D_d=\R_+ N_{d,w}$. The remaining conditions for a monoidal complex
are fulfilled as well. It is important to note that the monoidal
complex $\Mcc_{w}$ is not only dependent on $\Gamma_w$ or on the
pair $(\Gamma_w,E)$, but also on the chosen weight $w$.

The algebra $K[\Mcc_w]$ is again a residue class ring of the
polynomial ring $K[X_e:e\in E]$ under the as\-sign\-ment
$$
X_e \mapsto
\begin{cases}
t^{a_e} & \text{if $a_e\in |\Mcc_w|$},\\
0   & \text{else}.
\end{cases}
$$
The kernel of this epimorphism is denoted by $J_{\Mcc_{w}}$. It is
of course just the presentation ideal of the toric face ring
$K[\Mcc_w]$ supported by the conical complex $\Gamma_w$ (and here we
must allow that indeterminates $X^e$ go to $0$).


One cannot expect that $\ini_w(I_{\Mcc})=J_{\Mcc_w}$ since
$J_{\Mcc_w}$ is always a radical ideal, but $\ini_w(I_{\Mcc})$ need
not be radical. However, this is the only obstruction. The next
theorem generalizes a result of Sturmfels (see \cite{ST91} and
\cite{ST96}) who proved it in the case that conical complex is
induced from a single monoid and that the subdivision $\Gamma_w$ is
a triangulation. It is essentially equivalent to \cite[Theorem
5.11]{BRRO05}. See Remark \ref{othersources} for the difference of
the two approaches.

\begin{thm}
\label{Sturm1} Let $\Mcc$ be a monoidal complex supported on a
conical complex  $\Sigma$, let $(a_e)_{e\in E}$ be a family of
elements of $|\Mcc|$ such that $\{a_e:e\in E\}\cap M_c$ generates
$M_c$ for each $c\in\Sigma$, and let $w=(w_e)_{e\in E}$ be a weight
vector. Moreover, let $\Mcc_w$ be the monoidal complex defined by
$w$. Then the ideal $J_{\Mcc_w}$ is the radical of the initial ideal
$\ini_w(I_{\Mcc})$.
\end{thm}

\begin{proof}
For a single monoid the theorem is \cite[Theorem 7.18]{BRGU08}, and
we reduce the general case to it.

As remarked above, the ideal $J_{\Mcc_w}$ is the presentation ideal
of a toric face ring by construction. The underlying complex is
$\Gamma_w$, a subdivision of $\Sigma$. We apply Proposition
\ref{Mc_retr} to this subdivision of $\Sigma$ and the facets of
$\Sigma$. The latter correspond to single monoids $M_1,\dots,M_n$.
Thus
\begin{equation}
J_{\Mcc_w}=A_\Sigma+J_{(M_1)_w}+\dots+J_{(M_n)_w}.\label{Equa_J}
\end{equation}
By \cite[Theorem 7.18]{BRGU08} we have $J_{(M_1)_w}=\Rad
\ini_w(I_{M_i})$, and therefore
$$
J_{\Mcc_w}=A_\Sigma+\Rad S\cdot\ini_w(I_{M_1})+\dots+\Rad
S\cdot\ini_w(I_{M_n}).
$$
The right hand side is certainly contained in $\Rad
\ini_w(I_{\Mcc})$, and contains $\ini_w(I_{\Mcc})$ by Proposition
\ref{ini_retr}. Since $J_{\Mcc_w}$ is a radical ideal, we are done.
\end{proof}

Because of the equation $\Rad\ini_w(I_{\Mcc})=J_{\Mcc_w}$ we have
always the inclusion $\ini_w(I_{\Mcc}) \subseteq J_{\Mcc_w}$. It is
a natural question to characterize the cases in which we have
equality.  It holds exactly when the monoids $N_{d,w}$ are
determined by their cones:

\begin{cor}
\label{St_rad} With the hypotheses of Theorem \ref{Sturm1}, the
following statements are equivalent:
\begin{enumerate}
\item
$\ini_w(I_{\Mcc})$ is a radical ideal;
\item
For all facets $d\in \Gamma_w$ one has $N_{d,w}=M_c \cap D_d$ where
$c \in \Sigma$ is the smallest face such that $d\subseteq c$.
\end{enumerate}
\end{cor}

\begin{proof}
Condition (ii) evidently depends only on the facets of $\Sigma$, but
this holds for (i) likewise. The equality
$J_{\Mcc_w}=\ini_w(I_\Mcc)$ is passed to the facets, since we obtain
the corresponding ideals for the facets $c_i$ by intersection with
$S_{c_i}$, and in the converse direction we use equation
\eqref{Equa_J} and Proposition \ref{ini_retr}.
Therefore it is enough to consider the case of a single cone, in
which the corollary is part of \cite[Corollary 7.20]{BRGU08}.
\end{proof}

Before presenting another corollary we have to characterize the
cases in which $\ini_w(I_{\Mcc})$ is a monomial ideal. We say that a
monoidal complex is \emph{free} if all its monoids are free
commutative monoids. Evidently this implies that the associated
conical complex is simplicial, but not conversely. The free monoidal
complexes are exactly those derived from abstract simplicial
complexes (compare Example \ref{ExMonComp} (ii)). We note the
following obvious consequence of Theorem \ref{Sturm1}:

\begin{lem}
\label{St_rad2} Under the hypothesis of Theorem \ref{Sturm1} the
following statements are equivalent:
\begin{enumerate}
\item
$\Rad \ini_w(I_{\Mcc})$ is a (squarefree) monomial ideal;
\item
$\Mcc_w$ is a free monoidal complex.
\end{enumerate}
In particular, if these equivalent conditions hold, then $\Gamma_w$
is a regular triangulation of $\Sigma$.
\end{lem}

For the next result we recall the definition of unimodular cones.
Let $L\subseteq \R^d$ be a lattice, \ie  $L$ is a subgroup of $\R^d$
generated by $\R$-linearly independent elements, and we assume that
$L \subseteq \Q^d$. Let $C \subseteq \R^d$ be a rational pointed
cone. Since for each extreme ray $R$ of $C$ the monoid $R\cap L$ is
normal and of rank $1$, there exists a unique generator $e$ of this
monoid. We call these generators the {\em extreme generators} $C$
with respect to $L$. If $C$ is simplicial, then we call $C$ {\em
unimodular} with respect to $L$ if the sublattice of $L$ generated
by the extreme generators of $C$ with respect to $L$ generate a
direct summand of $L$.

\begin{thm}
\label{St_free} With the same assumptions as in Theorem \ref{Sturm1}
the following statements are equivalent:
\begin{enumerate}
\item
The ideal $\ini_w(I_{\Mcc})$ is a monomial radical ideal;
\item
The conical complex  $\Gamma_w$ is a triangulation of $\Sigma$, the
extreme generators of a cone $D_d$ for $d \in \Gamma_w$ with respect
to $\gp(M_c)$ generate the monoid $N_{d,w}$, and $D_d$ is unimodular
with respect to $\gp(M_c)$.
\end{enumerate}
\end{thm}
\begin{proof}
It follows from \ref{St_rad}  and \ref{St_rad2} that
$\ini_w(I_{\Mcc})$ is a monomial radical ideal if and only if:
\begin{enumerate}
\item[(a)]
$\Mcc_w$ is free.
\item[(b)]
For a facet $d\in \Gamma_w$ let $c \in \Sigma$ be the smallest face
such that $d\subseteq c$. Then $N_{d,w}=M_c \cap D_d$.
\end{enumerate}

It remains to show the equivalence of (a) and (b) to (ii). But both
sides of this equivalence depend only on the single monoids $M_c$
and the restrictions of $\Mcc_w$ to them. In the case of a single
monoid the theorem is part of \cite[Corollary 7.20]{BRGU08}.
\end{proof}

Now we can give a nice criterion for the normality of the monoids in
a monoidal complex in terms of an initial ideal with respect to a
weight vector.

\begin{thm}
\label{normal_nice} The following statements are equivalent:
\begin{enumerate}
\item
All monoids $M_c$ of the monoidal complex $\Mcc$ are normal.
\item
There exists a  family of elements $(a_e)_{e\in E}$  of $|\Mcc|$
such that $\{a_e: a_e\in M_c\}$ generates $M_c$ for each monoid
$M_c$ of $\Mcc$ and a weight vector $w=(w_e)_{e\in E}$ such that
$\ini_w(I_\Mcc)$ is a monomial radical ideal.
\end{enumerate}
\end{thm}
\begin{proof}
(ii) implies (i) is again reduced to the case of a single monoid $M$
by Proposition \ref{ini_retr}. In this case (ii) implies that $M$ is
the union of free monoids with the same group as $M$. Then the
normality of $M$ follows immediately.

For the converse we have to construct a regular unimodular
triangulation $\Gamma$ of $\R_+M$ by elements of $M$, which we
choose as a system of generators. The weight of 
$X_e$ is then chosen as the value of the support function of the
triangulation at $a_e$.

The existence of such a triangulation is a standard result. E.~g.\
see \cite[Theorem 2.70]{BRGU08} where it is stated for a single
monoid $M$. The construction goes through for monoidal complexes as
well (and the proof implicitly makes use of this fact). However,
there is one subtle point to be taken into account: if $M$ is normal
and $F$ is a face of the cone $\R_+M$, then $\gp(M\cap F)=\gp(M)\cap
\R F$. This condition ensures that the groups $\gp(M_c)$ form again
a monoidal complex, and that unimodularity of a free submonoid does
not depend on the monoid $M_c$ in which it is considered.
\end{proof}

In the investigation of a normal monoid $M$ one is usually not
interested in an arbitrary system of generators of $M$, but in
$\Hilb(M)$. It is well-known that one can not always find a
(regular) unimodular triangulation by elements of $\Hilb(M)$, and
this limits the value of results like Theorem \ref{normal_nice}
considerably. Nevertheless, it is very powerful when the
unimodularity of certain triangulations is given automatically.

Theorem \ref{normal_nice} can be used to prove that monoid algebras
of normal affine monoids are Cohen-Macaulay. This result is due to
Hochster \cite{HO72}.

\begin{cor}
Let M be a normal affine monoid. Then the monoid algebra $K[M]$
is Cohen-Macaulay for every field $K$.
\end{cor}

\begin{proof}
We may assume that $M$ is positive. In fact, $M=U(M)\oplus M'$ where
$U(M)$ is the group of units of $M$ and $M'$ is a normal affine
monoid which is positive. Moreover, $K[M]$ is a Laurent polynomial
extension of $K[M']$ and thus we may replace $M$ by $M'$.

It follows from Theorem \ref{normal_nice} that there exists a system
of generators $(a_e)_{e\in E}$  of $M$ and a weight vector
$w=(w_e)_{e\in E}$ such that $K[M]=S/I_{M}$ where $S=K[X_e :e \in
E]$ and $\ini_w(I_{M})$ is a monomial radical ideal. Thus
$\ini_w(I_{M})=I_\Delta$ for an abstract simplicial complex $\Delta$
on the vertex set $E$. Now standard results from Gröbner basis
theory yield that $K[M]$ is Cohen-Macaulay if the Stanley-Reisner
ring $K[\Delta]=S/I_\Delta$ is Cohen-Macaulay.

Observe that $\Delta$ is a triangulation of a cross section of $\R_+
M$. Now one can use \eg a theorem of Munkres \cite[5.4.6]{BRHE98}
which states that the Cohen-Macaulay property of $K[\Delta]$ only
depends on the topological type of $|\Delta|$. A cross-section of a
pointed cone is homeomorphic to a simplex whose Stanley-Reisner ring
is certainly Cohen-Macaulay.
\end{proof}

\section{Betti numbers of toric face rings}
\label{bettinumbers}

A consequence of Proposition \ref{presentation} is a presentation of
a toric face ring $K[\Mcc]$ over a polynomial ring $S$. It is a
natural question to determine the Betti numbers of $K[\Mcc]$ over
$S$, and the graded Betti numbers if there exists a natural grading.
The first question is of course which grading is a natural one to
consider. Recall that in general $|\Mcc|$ has only a partial monoid
structure and can not be used directly. But even if $\Sigma$ is a
fan in $\R^d$ and the monoids in $\Mcc$ are embedded in $\Z^d$, then
$\Z^d$ may not be the best choice to start with.

At first we recall a few facts from graded homological algebra. Let
$H$ be an (additive) commutative monoid which is positive, \ie $H$
has no invertible elements except $0$. Usually one defines graded
structures on rings and modules via groups. If $H$ is {\em
cancellative}, \ie  if $a+b = a+c$  implies  $b = c$ for $a,b,c \in
H$, then $H$ can be naturally embedded into the abelian
(Grothendieck) group $G$ of $H$. Therefore one can defines terms
like $H$-graded by considering $G$-graded objects whose homogeneous
components with degrees not in  $H$ are zero. But we will have to
consider noncancellative monoids, and thus it may be impossible to
embed $H$ into a group.

Hence we introduce $H$-graded objects directly. Let $R$ be a
commutative ring and $M$ an $R$-module (where we as always assume
that $R$ is commutative and not trivial). An {\em $H$-grading} of
$R$ is a decomposition $R=\bigoplus_{h \in H} R_h$ of $R$ as abelian
groups such that $R_h \cdot R_g \subseteq R_{h+g}$ for all $h,g \in
H$. A graded ring together with an $H$-grading is called an {\em
$H$-graded ring}. Now assume that $R$ is an $H$-graded ring. A {\em
grading} of $M$ is a decomposition $M=\bigoplus_{h \in H} M_h$ of
$M$ as abelian groups such that $R_h \cdot M_g \subseteq M_{h+g}$
for all $h,g \in H$. An $H$-graded $R$-module $M$ together with an
$H$-grading is called an {\em $H$-graded module}. $M_h$ is called
the {\em $h$-homogeneous component} of $M$ and an element $x \in
M_h$ is said to be {\em homogeneous} of degree $\deg x=h$.

From now on we assume that $R$ is a Noetherian $H$-graded ring. The
\fg $H$-graded $R$-modules build a category. The morphisms are the
homogeneous $R$-module homomorphisms $\phi \: M \to M'$, \ie
$\phi(M_h) \subseteq M'_h$ for all $h \in H$. For $h\in H$ we let
$M(-h)$ be the $H$-graded $R$-module with homogeneous components
$M(-h)_g=\bigoplus_{h' \in H, g=h'+h} M_{h'}$ for $g\in H$. In
particular, $R(-h)$ is a free $R$-module of rank 1 with
generator sitting in degree $h$. Since kernels of homogeneous maps
of \fg $H$-graded $R$-modules are again \fg $H$-graded and there
exist \fg free $H$-graded $R$-modules, every \fg $H$-graded
$R$-module has a free (hence projective) resolution
$$
F_{\lpnt} \: \dots \to F_n \to \dots \to F_0 \to 0
$$
where $F_n$ is a finite direct sum of free modules of the form
$R(-h)$ for some $h \in H$ and all maps are homogeneous and
$R$-linear.

Next we want to pose a condition on $H$ and specialize the
considered class of rings. We say that $H$ is {\em cancellative with
respect to $0$} if $a+b = a$  implies  $b = 0$ for $a,b \in H$. Let
$K$ be a field. An $H$-graded $K$-algebra $R$ is a Noetherian
$K$-algebra $R=\bigoplus_{h \in H}R_h$ with $R_0=K$. Since $H$ is
positive all homogeneous units of $R$ must belong to $R_0$ and $R$
has the unique $H$-graded maximal ideal $\mm=\bigoplus_{h \in
H\setminus\{0\}}R_h$. We see that $R$ is an $H$-graded local ring, a
notion defined in the obvious way. Observe that $\mm$ is also
maximal in $R$. The ring $R$ behaves like a local ring because of
the next lemma.

\begin{lem}
\label{nakayama} Assume that $H$ is cancellative with respect to $0$
and $R$ is an $H$-graded $K$-algebra. Then Nakayama's lemma hold,
\ie if $M$ is a \fg $H$-graded module and $N \subseteq M$ is a \fg
$H$-graded submodule s.t. $M=N+\mm M$, then $M=N$.
In particular,
homogeneous elements $x_1,\dots,x_n$ are a minimal system of generators of $M$
if and only if their residue classes
are a $K$-vector space basis of $M/\mm M$ and then we write $n=\mu(M_\mm)$.
\end{lem}

\begin{proof}
We may assume without loss of generality that $N=0$. Now
let $x_1,\dots,x_n$ be a minimal system of generators of homogeneous
elements of $M$. Since $M=\mm M$ we have an equation
$$
x_n= \sum_{i=1}^{n} a_i x_i
$$
where $a_i \in \mm$.
For all homogeneous components $a_{ij}$ of  some $a_i$ we may without loss of generality
assume that $\deg x_n=\deg a_{ij} +\deg x_i$. Fix a
homogeneous components $a_{in}$ of  $a_n$. Then
$\deg x_n=\deg a_{in} +\deg x_n$ implies $\deg a_{in}=0$ because $H$ is cancellative with
respect to $0$. Thus $a_{in} \in K\cap \mm$, and so $a_{in}=0$.
Hence $a_n=0$ and
$x_n$  is a linear combination of $x_1,\dots,x_{n-1}$, in
contradiction to the minimality of the system of generators.
\end{proof}

\begin{ex}
Let $F=\N^n$ for some $n\ge 0$, and set $\deg a=\sum a_i$ for $a\in
F$. Let $M$ a quotient of $F$ by a homogeneous congruence, \ie a
congruence in which $x\sim y$ implies $\deg x=\deg y$. Then $M$ is
cancellative at $0$, but in general it is not cancellative.
\end{ex}

Now we can re-prove many well-known results from local ring or
$\Z$-graded ring theory. (\Eg See \cite[Section 1.5]{BRHE98}.)

For example, let $x_1,\dots,x_n$ be a minimal system of generators
of $M$ where $\deg x_i=h_i\in H$, let $\phi \: F=\bigoplus_{i=1}^n
R(-h_i) \to M$ be the homogeneous map sending the generator $e_i$ of
$R(-h_i)$ to $x_i$.  Then we claim that $\Ker \phi \subseteq \mm F$.
Indeed, otherwise it follows that the residue classes of
$x_1,\dots,x_n$ are not a $K$-vector space basis of $M/\mm M$ and thus by
Nakayama's lemma $x_1,\dots,x_n$ is not minimal system of generators.
This is a contradiction.
Consequently there exist minimal $H$-graded
free resolutions, \ie  given a \fg $H$-graded module $M$, there
exists an $H$-graded free resolution $F_{\lpnt}$ of $M$ such that
$\Ker
\partial _n \subseteq \mm F_n$ for all $n$. If we write
$F_n=\bigoplus_{h \in H} R(-h)^{\beta^R_{n,h}(M)}$, then we call the
$\beta^R_{n,h}(M)$ the {\em $H$-graded Betti numbers} of $M$. Up to
homogeneous isomorphism  of complexes, $F_{\lpnt}$ is uniquely
determined by the requirement that $\Ker \partial _n \subseteq \mm
F_n$ all $n$. The numbers $\beta^R_{n,h}(M)$ are also uniquely
determined. Indeed, $\Tor^R_n(M,K)$ is an $H$-graded module
considered as an $R$- or $K=R/\mm$-module and we have $\dim_K
\Tor^R_n(M,K)_h=\beta^R_{n,h}(M)$ as is easily verified. Some more
results in this direction can easily be verified: a \fg $H$-graded
$R$-module is projective if and only if it is free, $\projdim
M=\projdim M_\mm$ and so forth.

Next we want to apply the theory discussed so far, to the situation
of toric face rings. Let $\Mcc$ be a monoidal complex supported on a
pointed conical complex $\Sigma$, and let $(a_e)_{e\in E}$ be a
family of elements of $|\Mcc|$ such that $\{a_e:e\in E\}\cap M_c$
generates $M_c$ for each $c\in\Sigma$.

According to Proposition \ref{presentation} the defining ideal
$I_\Mcc$ of the toric face ring $K[\Mcc]$ of a monoidal complex is a
sum
$$
I_\Mcc=A_\Mcc+B_\Mcc
$$
where $A_\Mcc$ is an ideal generated by squarefree monomials and
$B_\Mcc$ is a binomial ideal containing no monomials: This is a
consequence of the fact that every binomial generator vanishes on
the vector $(1)_{e \in E}$, but a monomial has here value $1$.

Recall that a congruence relation on a commutative monoid $M$ is an
equivalence relation $\sim$ such that for $a,b,c \in M$ with $a\sim
b$ we have $a+c\sim b+c$. Now $M/\sim $ is again a commutative
monoid in a natural way.

Consider the free monoid $\N^E$ with generators $f_e$ for $e \in E$.
Note that $S=K[X_e : e \in E]$ is the monoid algebra of $\N^E$; the
monomials in $S$ are denoted by $X^a=\prod_{e\in E}X_e^{a_e}$. On
$\N^E$ we define the congruence relation $a \sim b$ for $a,b \in
\N^E$ if and only if $X^a-X^b \in B_\Mcc$ is a binomial.
We let $H_{\Mcc}$ denote the monoid $\N^E/\sim$. It is
not to hard to see and well-known that $S/B_\Mcc$ is exactly the
monoid algebra of the monoid $H_{\Mcc}$.

\begin{lem}
\label{positivemonoid} $H_{\Mcc}$ is a commutative positive monoid
with monoid algebra $S/B_\Mcc$.
\end{lem}

\begin{proof}
It only remains to show that $H_{\Mcc}$ is positive. Let
$\overline{g},\overline{h} \in H_{\Mcc}$ for $g,h \in \N^E$ such
that $\overline{g+h}=\overline{0}$ and assume that
$\overline{g},\overline{h} \neq \overline{0}$. It follows from the
definition of $H_{\Mcc}$ that $X^{g+h}-1 \in B_\Mcc$. But $B_\Mcc$
is generated by binomials that vanish on the zero vector $(0)_{e \in
E}$ because all monoids $M_c$ for $c \in \Sigma$ are positive. The
binomial $X^{g+h}-1$ does not vanish on $(0)$, and this yields a
contradiction.
\end{proof}

We saw that from the algebraic point of view it is very useful if
$H_{\Mcc}$ is cancellative with respect to $0$. But it is not strong
enough for a combinatorial description of the Betti numbers, as a
counterexample will show. The next lemma describes a stronger
cancelation property for monoidal complexes associated with fans.

\begin{lem}
\label{cancel} Assume that $\Sigma$ is a rational pointed fan  in
$\R^n$ such that $M_c \subseteq \Z^n$ for $c \in \Sigma$.
\begin{enumerate}
\item
If $\overline{i}+\overline{j} =\overline{i}+\overline{k}$ for
$\overline{i},\overline{j}, \overline{k} \in H_{\Mcc}$ where $i,j,k
\in \N^E$, then $X^j- X^k \in I_{\Mcc}$.
\item
The monoid $H_{\Mcc}$ is cancellative with respect to $0$.
\item
If $X^j- X^k \in I_{\Mcc}$ and $X^j, X^k \not\in I_{\Mcc}$, then
$\overline{i}=\overline{j}$ in $H_{\Mcc}$.
\end{enumerate}
\end{lem}

\begin{proof}
(i) Note that the toric face ring $K[\Mcc]$ has a natural
$\Z^n$-grading induced by the embeddings $M_c \subseteq \Z^n$ for $c
\in \Sigma$. Then also the polynomial ring $K[X_e: e\in E]$ is
$\Z^n$-graded if we give $X_e$ the degree $a_e \in \Z^n$. Observe
that the ideal $B_\Mcc$ is then $\Z^n$-graded since the generators
are homogeneous with respect to this grading. Then $K[X_e: e\in
E]/B_\Mcc$ is $\Z^n$-graded. Equivalently we obtain a monoid
homomorphism $\phi\: H_{\Mcc} \to \Z^n, \overline{i} \mapsto \sum_{e
\in E}i_e a_e$.

Now $\overline{i}+\overline{j}=\overline{i}+\overline{k}$ implies
that $\sum_{e \in E} (i_e +j_e)a_e=\sum_{e \in E} (i_e+k_e)a_e$ in
$\Z^n$, and thus $\sum_{e \in E} j_e a_e = \sum_{e \in E} k_e a_e$.
It follows that $X^j- X^k \in \Ker (K[X_e : e\in E] \to K[\Mcc])=
I_{\Mcc}$.

(ii) It follows from (i) that $H_{\Mcc}$ is cancellative with
respect to $0$, because $X^j- 1 \not\in I_{\Mcc}$.

(iii) If $X^j- X^k \in I_{\Mcc}$ and $X^j, X^k \not\in I_{\Mcc}$,
then  $X^j- X^k \in B_\Mcc$ by the last observation of Proposition
\ref{presentation}. Hence $\overline{i}=\overline{j}$ in $H_{\Mcc}$.
\end{proof}

Let $S=K[X_e : e\in E]$. Observe that all rings $S, S/B_\Mcc$ and
$K[\Mcc]$ are naturally $H_{\Mcc}$-graded. For $S$ we set $\deg
X^i=\overline{i}\in H_{\Mcc}$. Since $X^a-1 \not\in I_{\Mcc}$ we
have that $H_\Mcc$ is positive, and $S$ is an $H$-graded local ring.
If $H_{\Mcc}$ is cancellative with respect to $0$, then one can
apply Lemma \ref{nakayama} to \fg $H_{\Mcc}$-graded $S$-modules like
$K[\Mcc]$. In particular, we can speak about minimal
$H_{\Mcc}$-graded resolutions. The next goal is to determine the
corresponding $H_{\Mcc}$-graded Betti numbers.

For $\overline{h} \in H_{\Mcc}$ we define
$$
\Delta_{\overline{h}} = \{ F\subseteq E : \overline{h} =
\overline{g} + \sum_{e \in F} \overline{f_e} \text{ for some }
\overline{g} \in  H_{\Mcc}  \}.
$$
We see immediately that $\Delta_{\overline{h}}$ is a simplicial
complex on the finite vertex set $E$, and we call
$\Delta_{\overline{h}}$ the {\em squarefree divisor complex} of
$\overline{h}$. Moreover, we need a special subcomplex of
$\Delta_{\overline{h}}$ defined as follows
$$
\Delta_{\overline{h}, \Mcc} = \{ F\subseteq E : \overline{h} =
\overline{g} + \sum_{e \in F} \overline{f_e} \text{ for some }
\overline{g} \in  H_{\Mcc} \text{ such that } X^g \in I_{\Mcc}\}.
$$
For an arbitrary simplicial complex $\Delta$ on some ordered vertex
set $E$ (with order $<$) we let  $\Tilde \Cc(\Delta)$ denote the
augmented oriented chain complex of $\Delta$ with coefficients in
$K$, \ie  the complex
$$
\Tilde \Cc_{\lpnt}(\Delta) \: 0\to \Cc_{\dim \Delta}
\overset{\partial}{\to} \dots \overset{\partial}{\to} \Cc_0
\overset{\partial}{\to} \Cc_{-1} \to 0
$$
where  $\Cc_{i} = \bigoplus_{F \in \Delta,\ \dim F=i} KF $ and
$\partial (F)= \sum_{F' \in \Delta,\ \dim F'=i-1} \epsilon(F,F')
F'$. Here $\epsilon(F,F')$ is $0$ if $F' \not\subseteq F$. Otherwise
it is $(-1)^{k}$ if $F=\{ e_0, \dots, e_{i} \}$ for elements $e_0<
\dots< e_{i}$ in $E$, and $F'=\{ e_0, \dots,e_{k-1},e_{k+1},\dots,
e_{i} \}$. Further we let $\Tilde H(\Delta)_{i}=H_i(\Tilde
\Cc(\Delta)_{\lpnt})$ be the $i$-th reduced simplicial homology
group of $\Delta$. If $\Delta'$ is a subcomplex of $\Delta$ we let
$\Tilde \Cc(\Delta,\Delta')_{\lpnt}=\Tilde
\Cc(\Delta)_{\lpnt}/\Tilde \Cc(\Delta')_{\lpnt}$ denote the relative
augmented oriented chain complex of $\Delta$ and $\Delta'$, and
$H_{i}(\Delta,\Delta')$ the $i$-th homology of this complex.

\begin{thm}
\label{tor} Assume that $\Sigma$ is a rational pointed fan in $\R^n$
such that $M_c \subseteq \Z^n$ for $c \in \Sigma$. Let $(a_e)_{e\in
E}$ be a family of elements of $|\Mcc|$ such that $\{a_e:e\in
E\}\cap M_c$ generates $M_c$ for each $c\in\Sigma$. Then
$$
 \beta_{ih}^S(K[\Mcc]) =
 \dim_K \Tilde H_{i-1} (\Delta_{\overline{h}},\Delta_{\overline{h}, \Mcc})
$$
for all $\overline{h} \in H_{\Mcc}$ and $i \in \N$.
\end{thm}

\begin{proof}
Let $S=K[X_e : e \in E]$. We fix an arbitrary total order $<$ on
$E$. Let $\Kc_{\lpnt}(K[\Mcc])$ denote the Koszul complex of $X_e,
e\in E$ tensored with $K[\Mcc]$. This complex is naturally
$H_{\Mcc}$-graded and its $H_{\Mcc}$-graded homology is exactly
$\Tor_{\lpnt}(K,K[\Mcc])$. (\Eg see \cite{BRHE98} for details.) Thus
we can use this complex to determine the numbers
$\beta_{ih}^S(K[\Mcc])$. We have
$$
\Kc_{i}(K[\Mcc]) = \bigoplus_{F \subseteq E,\ |F|=i} K[\Mcc]
(-\sum_{e\in F} \overline{f_e})
$$
and the differential $\partial_i \: \Kc_{i}(K[\Mcc]) \to
\Kc_{i-1}(K[\Mcc])$ is given on the component
$$
K[\Mcc] (-\sum_{e\in F} \overline{f_e}) \to K[\Mcc] (-\sum_{e\in F'}
\overline{f_e})
$$
for $F', F \subseteq E$ as the zero map for $F'\not\subseteq F$, or
otherwise as multiplication $\epsilon(F,F')X_{e_k}$ where
$$
\epsilon(F,F')
=
\begin{cases}
0 & \text{if } F' \not\subseteq F,\\
(-1)^{k-1} & \text{if } F=\{e_1<\dots<e_i\},  F'=F \setminus \{e_k\}.
\end{cases}
$$
For $\beta_{i\overline{h}}^S(K[\Mcc])$ we have first to determine
$\Kc_{i}(K[\Mcc])_{\overline{h}}$. Thus we compute
$$
K[\Mcc] (-\sum_{e\in F} \overline{f_e})_{\overline{h}} =
\bigoplus_{\overline{h'} \in H_{\Mcc},\ \overline{h'}+\sum_{e\in F}
\overline{f_e}= \overline{h}} K[\Mcc]_{\overline{h'}}.
$$
Such an $\overline{h'}$ exists if and only if $F \in
\Delta_{\overline{h}}$. If $K[\Mcc]_{\overline{h'}} \neq 0$, then
$X^{h'} \not\in I_{\Mcc}$. Assume that there exists another
$\overline{h''}$ such that $\overline{h''}+\sum_{e\in F}
\overline{f_e}= \overline{h}$ and $K[\Mcc]_{\overline{h'}} \neq 0$.
It follows $X^{h''} \not\in I_{\Mcc}$. We obtain from \ref{cancel}
and $\overline{h'}+\sum_{e\in F}
\overline{f_e}=\overline{h''}+\sum_{e\in F} \overline{f_e}$ that
$X^{h'}-X^{h''} \in B_\Mcc$. Hence $\overline{h'}=\overline{h''}$ in
$H_{\Mcc}$, \ie $\overline{h'}$ is uniquely determined if it exists.
Moreover, we have $K[\Mcc]_{\overline{h'}} =0$ if and only if $F \in
\Delta_{\overline{h}, \Mcc}$. Hence
$$
\Kc_{i}(K[\Mcc]) \cong \bigoplus_{F \in \Delta_{\overline{h}}
\setminus \Delta_{\overline{h}, \Mcc},\ |F|=i} K F.
$$
Consider $F' \in \Delta_{\overline{h}} \setminus
\Delta_{\overline{h}, \Mcc}$ such that $|F'|=i-1$. Choose
$\overline{h'}$ such that $\overline{h'}+\sum_{e\in F}
\overline{f_e}= \overline{h}$ and $\overline{h''}$ such that
$\overline{h''}+\sum_{e\in F'} \overline{f_e}= \overline{h}$. The
differential $\Kc_{i}(K[\Mcc]) \to \Kc_{i-1}(K[\Mcc])$ on the
component $KF\to KF'$ (which corresponds to $K[\Mcc]_{\overline{h'}}
\to K[\Mcc]_{\overline{h''}}$) is given by
$$
\partial_i(F)=
\begin{cases}
0 & \text{if } F' \not\subseteq F,\\
\epsilon(F,F')F'& \text{if } F=\{e_1<\dots<e_i\},  F'=F \setminus \{e_k\}.
\end{cases}
$$
Hence we see that the complex $\Kc_{\lpnt}(K[\Mcc])_{\overline{h}}$
coincides with $\Tilde \Cc_{\lpnt-1}(
\Delta_{\overline{h}},\Delta_{\overline{h}, \Mcc})$ and this yields
$$
 \beta_{ih}^S(K[\Mcc]) =
 \dim_K \Tilde H_{i-1} (\Delta_{\overline{h}},\Delta_{\overline{h}, \Mcc})
$$
as desired.
\end{proof}

One can easily generalize Theorem \ref{tor} in the following way: if
$\Mcc$ satisfies the properties of Lemma \ref{cancel}, then the
proof of Theorem \ref{tor} works for $\Mcc$. However, in general one
can not expect that the compact combinatorial formula is true for
all monoidal complexes without any further assumptions. Indeed, a
counterexample is

\begin{ex}
\label{Moeb_mon} We consider the Möbius strip as a monoidal complex
$\Mcc$ by considering each quadrangle as a unit square and choosing
the monoid over it as the corresponding monoid. Together with the
compatibility conditions this determines $\Mcc$ completely.
\begin{figure}[hbt]
$$
 \psset{unit=2cm}
\begin{pspicture}(-0.7,-0.5)(1.3,1)
\def\vertex{\pscircle[fillstyle=solid,fillcolor=black]{0.05}}
 \psset{viewpoint=3.5 1 1.5 linewidth=1.5pt}
\ThreeDput[normal=0 0 1](0,0,0){
 \pspolygon(0,0)(0,1)(1,1)(1,0)
 \rput(0,0){\vertex}
 \rput(0,1){\vertex}
 \rput(1,0){\vertex}
 \rput(1,1){\vertex}
 }
\ThreeDput[normal=0 -1 0](0,0,0){
 \pspolygon(0,0)(0,1)(1,1)(1,0)
 \rput(0,1){\vertex}
 \rput(1,1){\vertex}
 }
\ThreeDput[normal=0 1 1](1,1,0){
 \pscurve(0,0)(0.1,0.4)(0.43,0.62)
 \pscurve(0.64,0.68)(0.9,0.9)(1,1.414)
 \pscurve(0,1.414)(0.2,0.9)(0.9,0.4)(1,0)
 }
 \rput(-0.4,0.5){$x$}
 \rput(-0.4,-0.4){$y$}
 \rput(-0.2,0.9){$u$}
 \rput(-0.15,0.0){$v$}
 \rput(0.85,-0.45){$z$}
 \rput(1.15,-0.1){$w$}
\end{pspicture}
$$
\caption{M\"obius strip as a monoidal complex}
\label{FigMoeb_mon}
\end{figure}
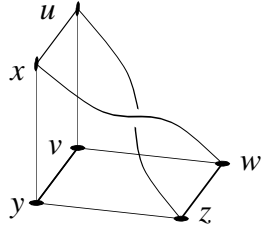

The ideal $I_{\Mcc}$ is generated by the binomials resulting from
the unit squares and monomials
$$
X_xX_z-X_uX_w,\quad X_yX_w-X_vX_z,\quad X_xX_v-X_uX_y,\quad X_u X_v
X_w \text{ and } X_u X_v X_z.
$$
The other monomials are redundant. \Eg $ X_u X_y X_z =X_z(X_u
X_y-X_x X_v) + X_v(X_x X_z - X_u X_w) + X_v X_u X_w$. Since the
binomial relations are homogeneous, $H_\Mcc$ is cancellative with
respect to $0$.

Let $x_a$ stand for the residue class of $X_a$, and choose the
degree
$$
\overline{h} = x_u x_v x_z = x_u x_w x_y = x_x x_y x_z = x_v x_w x_x
\in H_{\Mcc}.
$$
This equation shows that Lemma \ref{cancel} does not hold for
$\Mcc$.

Since $X^h\in I_\Mcc$, one has $K[\Mcc]_{\overline{h}}=0$. The
degree $\overline h$ component of the Koszul complex is
$$
\Kc_{\lpnt}(K[\Mcc])_{\overline{h}} \:
0 
\to
K^4 
\to
K^{12} 
\to
K^9 
\to
0 
\to 0
$$
where $K^9$ is in homological degree $1$, and the Betti numbers are
$$
\beta_{i \overline{h}}^S(K[\Mcc]) =
\begin{cases}
1 &\text{for } i=1,\\
0 &\text{else}.
\end{cases}
$$
Now we consider the complex $\Tilde
\Cc_{\lpnt}(\Delta_{\overline{h}},\Delta_{\overline{h},\Mcc})$ which
is given by
$$
0 
\to
K^4 
\to
K^{12} 
\to
K^6 
\to
0 
\to 0
$$
with $K^6$ in homological degree $0$. Hence $\Tilde
H_{1}(\Delta_{\overline{h}},\Delta_{\overline{h},\Mcc})=K^d$ for
some $d\geq 2$, and the formula of Theorem \ref{tor} dos not hold in
this case.
\end{ex}

The results of this section imply, in particular, the well-known Tor
formula of Hochster for Stanley-Reisner rings (see \cite{HO77}).
Recall that if $\Delta$ is a simplicial complex on the vertex set
$[n]$ and we let $S=K[X_1,\dots,X_n]$, then $K[\Delta]=S/I_\Delta$
is the Stanley-Reisner ring of $\Delta$
where $I_\Delta=(\prod_{i \in F}X_i : F \not\in \Delta)$
is the Stanley-Reisner ideal of $\Delta$.
Now all considered rings
have a natural $\Z^n$-grading. It is well-known that
$\beta_{ia}^S(K[\Delta])=0$ if $a$ is not a squarefree vector, \ie a
0-1 vector. (Either one shows this by using the results of this
section, or proves this directly). For a squarefree vector $a$ with
support $W=\{i\in[n]: a_i=1\}$, we write
$\beta_{iW}^S(K[\Delta])=\beta_{ia}^S(K[\Delta])$ for the
corresponding Betti-number.

\begin{cor}[Hochster]
Let $\Delta$ be a simplicial complex on the vertex set $[n]$. Then
for $W\subseteq [n]$ one has
$$
\beta_{iW}^S(K[\Delta])
=
\dim_K \Tilde H_{|W|-i-1} (\Delta_W)
$$
where $\Delta_W=\{F \in \Delta : F \subseteq W\}$.
\end{cor}

\begin{proof}
In Example \ref{ExMonComp} it was observed that there exists a
rational pointed fan $\Sigma$ and an embedded monoidal complex
$\Mcc$ such that $K[\Mcc]=K[\Delta]$. Thus the binomial ideal
$B_\Mcc$ is $0$, and $I_{\Mcc}=I_\Delta$ is generated by squarefree
monomials. The monoid $H_{\Mcc}$ is nothing but the free monoid
$\N^n$ in this case. Thus the induced grading is just the natural
$\N^n$-grading on $K[\Delta]$. It remains to observe that the
complex $\Tilde \Cc_{\lpnt-1}
(\Delta_{\overline{h}},\Delta_{\overline{h}, \Mcc})$ coincide with
the complex $\Tilde \Cc_{|W|-\lpnt-1}( \Delta_{W})$ which determines
the homology $\Tilde H_{|W|-i-1} (\Delta_W)$. This concludes the
proof.
\end{proof}

\begin{rem}
Let $\Delta$ be a simplicial complex on  $[n]$. Hochster computed
also the local cohomology of the Stanley-Reisner ring as a
$\Z^n$-graded $K$-vector space in terms of combinatorial data of the
given complex (for example, see \cite{BRHE98} or \cite{ST05}). While
the Tor formula is restricted to embedded monoidal complexes (or
complexes which behave like these), one can prove a Hochster formula
for the local cohomology in great generality. In fact, to show such
a formula for toric face rings was one of the starting points of the
systematic study of toric face rings. See \cite{ICRO} for the case
of embedded monoidal complexes and \cite{BRBRRO07} for classes of
rings which include toric face rings as a special case.
\end{rem}

\end{document}